\documentclass[11pt,letterpaper,twoside]{article}
\usepackage{amsmath,amssymb,amsthm,amsfonts,amstext,amsbsy,amscd}

\usepackage{comment}
\usepackage{graphicx}
\usepackage{mathrsfs}
\usepackage{harvard}

\input{mymath.sty}

\providecommand{\MR}{$\clubsuit$}

\allowdisplaybreaks[1]

\renewcommand{\cal}{\mathscr}     

\begin{document}
\title{Asymptotic equivalence for nonparametric regression with multivariate and random design}
\author{\parbox[t]{6.5cm}{\centering Markus
 Rei\ss\\[2mm]
 \normalsize{\it
Institute of Applied Mathematics\\
University of Heidelberg}\\
\mbox{} reiss@statlab.uni-heidelberg.de}} \maketitle

\begin{abstract}
We show that nonparametric regression is asymptotically equivalent
in Le Cam's sense with a sequence of Gaussian white noise
experiments as the number of observations tends to infinity. We
propose a general constructive framework based on approximation
spaces, which permits to achieve asymptotic equivalence even in the
cases of multivariate and random design.
\end{abstract}

{\small \noindent {\it Key words and Phrases:} Le Cam deficiency,
equivalence of experiments, approximation space, interpolation,
Gaussian white noise.

\newpage

\section{Introduction}

Nonparametric regression is the model most often encountered in
nonparametric statistics because of its widespread applications. On
the other hand, for theoretical investigations the Gaussian white
noise or sequence space model is often preferred since it exhibits
nice mathematical properties. The common wisdom that statistical
decisions in the two models show the same asymptotic behaviour has
been formalized and proved for the first time by \cit{BrownLow} in
the one-dimensional case, using Le Cam's concept of equivalence of
statistical experiments.

In this paper we propose a unifying framework for establishing
global asymptotic equivalence between Gaussian nonparametric
regression and white noise experiments based on constructive
transitions with only minimal randomisations. This framework not
only allows to give concise proofs of known results, but extends the
asymptotic equivalence to the multivariate and random design
situation. The multivariate result has often been alluded to,
though it has never been proved, see e.g. \cit{HoffLep}. 
While \cit{BrownZhang} remark that the regression and white noise
experiments are not asymptotically equivalent for equidistant design
on $[0,1]^d$ and Sobolev classes of regularity $s\le d/2$, the so
far only positive result by \cit{Carter} states asymptotic
equivalence for equidistant design in dimensions $d=2$ and $d=3$
when $s>d/2$. The difficulty in extending results to higher
dimensions is that we have to go beyond piecewise constant or linear
approximations. For the dynamic model of ergodic diffusions
\cit{DalReiss2} have established multidimensional asymptotic
equivalence with a white noise model. For the case of univariate
nonparametric regression, but with non-Gaussian errors we refer to
\cit{GramaNussbaum}.

In Section \ref{SecIsometric} the concept of isometric approximation
spaces is introduced and applied to local constant and Fourier
approximations. The latter yields an easy proof for asymptotic
equivalence in any dimension $d$ for periodic Sobolev classes of
regularity $s>d/2$ and extends scalar results by \cit{Rohde}. A more
flexible framework is obtained using isomorphic approximations
spaces in Section \ref{SecIsomorph}. As a main application, a
constructive asymptotic equivalence result is established on the
basis of wavelet multiresolution analyses, which provides
equivalence results also for non-periodic function classes.
Connections to asymptotic studies by \cit{DonJoh2} and
\cit{JohnSilv} for wavelet estimators are discussed. The case of a
random design, uniform on a $d$-dimensional cube, is treated in
Section \ref{SecRanDes}. This setting is much more involved, but can
also be cast in the isomorphic framework. The construction is based
on a two-level procedure, generalizing an idea by \cit{Brownetal}.
Fine approximation and symmetry properties of the Fourier basis
yield the main result that also in the case of random design
asymptotic equivalence holds for Sobolev regularities $s>d/2$ and
any dimension $d\ge 1$.

\section{Isometric approximation}\label{SecIsometric}

\subsection{General theory}

We write ${\cal L}^2({\cal D}):=\{f:{\cal D}\to
\K\,|\,\norm{f}_{L^2}^2:=\int \abs{f}^2<\infty\}$ with $\K=\R$ or
$\K=\C$ and $L^2({\cal D})$ for the Hilbert space of equivalence
classes with respect to $\norm{\cdot}_{L^2}$. Although the
observations are real-valued, we shall use complex-valued functions
for simplicity when treating Fourier approximations.

\begin{definition}
Let $\E^d_n$ be the regression experiment obtained from observing
\[ Y_i=f({x}_i)+\sigma\eps_i,\qquad i=1,\ldots,n,
\]
for $n\in\N$, $f:{\cal D}\subset \R^d\to\R$ in some class ${\cal
F}^d\subset {\cal L}^2({\cal D})$, fixed design points $x_i\in {\cal
D}$ and for independent random variables $\eps_i\sim {\cal N}(0,1)$.
\end{definition}

Suppose we are given an $n$-dimensional space $S_n\subset {\cal
L}^2({\cal D})$ and a linear mapping $D_n:{\cal L}^2({\cal
D})\to\K^n$ with the following isometric property on $S_n$
\begin{equation}\label{EqIsometric}
\forall g_n\in S_n:\, \norm{g_n}_{L^2}=\norm{g_n}_n
:=n^{-1/2}\abs{D_ng_n}_{\K^n}.
\end{equation}
By $\scapro{\cdot}{\cdot}_n$ we denote the scalar product associated
with $\norm{\cdot}_n$. Usually, $D_ng=(g(x_i))_{1\le i\le n}$ will
be the point evaluation at the $n$ design points in which case
$\norm{g}_n^2=\frac1n\sum_{i=1}^n\abs{g(x_i)}^2$ is just the
empirical norm. Let us further introduce the linear operator
\[ {\cal I }_n:{\cal L}^2({\cal D})\to S_n,\quad
{\cal I }_ng:=(D_n|_{S_n})^{-1}(g(x_1),\ldots,g(x_n))^\top.
\]
For $D_ng=(g(x_i))_{1\le i\le n}$ we have ${\cal
I}_n=(D_n|_{S_n})^{-1}D_n$ and ${\cal I}_n$ is the
$\norm{\cdot}_n$-orthogonal projection onto $S_n$ such that ${\cal
I}_ng$ is the unique element of $S_n$ interpolating $g$ at the
design points $(x_i)$.

To state the first results, we refer to \cit{LeCamYang} for the
notion of equivalence between experiments and of the Le Cam distance
between two experiments $\E$ and $\G$, which for the parameter class
$\cal F$ will be denoted by $\Delta_{\cal F}(\E,\G)$. The Gaussian
law on a Hilbert space $H$ with mean vector $\mu\in H$ and
covariance operator $\Sigma:H\to H$ will be denoted by ${\cal
N}(\mu,\Sigma)$.

The regression experiment $\E^d_n$ can be transformed to a
functional Gaussian shift experiment by applying the isometry
$(D_n|_{S_n})^{-1}$ to $Y=(Y_i)\in\R^n$:
\begin{equation}\label{EqZ}
Z:=(D_n|_{S_n})^{-1}Y={\cal I}_nf+\frac{\sigma}{\sqrt{n}}\zeta\;\in
S_n,
\end{equation}
where $\zeta:=\sqrt{n}(D_n|_{S_n})^{-1}\eps\sim{\cal
N}(0,\Id_{S_n})$ is a Gaussian white noise in $S_n$ because for
$g_n,h_n\in S_n$
\[ E[\scapro{\zeta}{g_n}_{L^2}\overline{\scapro{\zeta}{h_n}}_{L^2}]
=n^{-1}
E[\scapro{\eps}{D_ng_n}_{\K^n}\overline{\scapro{\eps}{D_nh_n}}_{\K^n}]
= \scapro{g_n}{h_n}_{n}=\scapro{g_n}{h_n}_{L^2}.
\]

By adding completely uninformative observations on the orthogonal
complement of $S_n$ in $L^2({\cal D})$, the observation of $Z$ in
\eqref{EqZ} is equivalent to observing
\[ \scapro{\phi}{Z}_{L^2}=\scapro{\phi}{{\cal
I}_nf}_{L^2}+\frac{\sigma}{\sqrt{n}}\scapro{\phi}{\bar\zeta}_{L^2}
\quad\forall\,\phi\in L^2({\cal D}),
\]
with $\scapro{\phi}{\bar\zeta}_{L^2}\sim {\cal
N}(0,\norm{\phi}_{L^2}^2)$. In differential notation we have thus
established the following equivalence.

\begin{proposition}\label{PropEquiv1}
Let $\F^d_n$ be the Gaussian white noise experiment in $L^2({\cal
D})$ given by observing
\[ dY(x)={\cal I }_nf(x)\,dx+\frac{\sigma}{\sqrt{n}}dB(x),\quad x\in {\cal D},
\]
where $f\in{\cal F}^d$ and $dB$ is a Gaussian white noise in
$L^2({\cal D})$. Then the regression experiment $\E^d_n$ is
statistically equivalent to $\F^d_n$ for any functional class ${\cal
F}^d$.
\end{proposition}

We are coming to the first main result.

\begin{definition}
Let $\G^d_n$ be the Gaussian white noise experiment given by
observing
\[ dY(x)=f(x)\,dx+\frac{\sigma}{\sqrt{n}}dB(x),\quad x\in {\cal D},
\]
where $f\in{\cal F}^d$ and $dB$ is a Gaussian white noise in
$L^2({\cal D})$.
\end{definition}

\begin{theorem}\label{ThmMain}
The Le Cam distance between $\E^d_n$ and $\G^d_n$ for the class
${\cal F}^d$ is bounded by
\[ \Delta_{{\cal F}^d}(\E^d_n,\G^d_n)\le 1-2\Phi\Big(-\frac{\sqrt{n}}{2\sigma}\sup_{f\in{\cal F}^d}
\norm{f-{\cal I }_nf}_{L^2}\Big),
\]
where $\Phi$ denotes the standard Gaussian cumulative distribution
function.
\end{theorem}

\begin{remark}\label{Pn}
Note that $\norm{f-{\cal
I}_nf}_{L^2}^2=\norm{f-P_nf}_{L^2}^2+\norm{P_nf-{\cal
I}_nf}_{L^2}^2$ holds where $P_n$ is the $L^2$-orthogonal projection
onto $S_n$. This means that the bound on the Le Cam distance is
always larger than the same expression involving the classical bias
estimate $\sup_{f\in{\cal F}^d}\norm{f-P_nf}_{L^2}$. Because of
$\Phi(0)=1/2$ Proposition \ref{ThmMain} yields the rate estimate
\[\Delta_{{\cal F}^d}(\E^d_n,\G^d_n)\lesssim \sigma^{-1}n^{1/2}\sup_{f\in{\cal F}^d}
\norm{f-{\cal I }_nf}_{L^2}.
\]
Here and in the sequel $A\lesssim B$ means $A\le c B$ with a
constant $c>0$, independent of the other parameters involved, and
$A\thicksim B$ is short for $A\lesssim B$ and $B\lesssim A$.
\end{remark}

\begin{proof}
Since $\E^d_n$ and $\F^d_n$ are equivalent, it suffices to establish
the bound for $\Delta_{{\cal F}^d}(\F^d_n,\G^d_n)$. The two latter
experiments are realized on the same sample space. Therefore the Le
Cam distance is bounded by the maximal total variation distance over
the class ${\cal F}^d$ \cite[Prop. 2.2]{Nussbaum}. For Gaussian
white noise the total variation distance is given by
$1-2\Phi(-\frac{\sqrt{n}}{2\sigma} \norm{f-{\cal I}_nf}_{L^2})$
\cite[Section 3.2]{Carter}, and the result follows.
\end{proof}

\subsection{Piecewise constant approximation}\label{SecConstAprrox}

The original results of \cit{BrownLow} for equidistant design on
${\cal D}=(0,1]$ fit into the proposed isometric framework. For
design points $x_i=i/n$, $i=1,\ldots,n$, we consider the
$n$-dimensional space $S_n$ of piecewise constant, left-continuous
functions on $(0,1]$ with possible jumps at $i/n$, $i=1,\ldots,n-1$.
Using $D_ng=(g(i/n))_{1\le i\le n}$ we obtain for $g_n\in S_n$
\[ \norm{g_n}_n^2=\frac1n\sum_{i=1}^n \abs{g_n(i/n)}^2=\sum_{i=1}^n\int_{(i-1)/n}^n
\abs{g_n(u)}^2\,du=\norm{g_n}_{L^2}^2
\]
such that $D_n$ has the isometric property. To infer asymptotic
equivalence by Proposition \ref{ThmMain}, we have to ensure that
$\norm{f-{\cal I}_nf}_{L^2}=o(n^{-1/2})$ uniformly over all $f$ in
some functional class ${\cal F}^d$. Considering the H\"older class
of regularity $\alpha\in (0,1]$
\[ {\cal F}_H(\alpha,R):=\Big\{ f\in C^\alpha([0,1])\,\Big|\, \sup_{x\not=
y}\abs{f(x)-f(y)}/\abs{x-y}^\alpha\le R\Big\},
\]
we obtain for $f\in {\cal F}_H(\alpha,R)$
\begin{align*}
\norm{f-{\cal I}_nf}_{L^2}^2& = \sum_{i=1}^n \int_{(i-1)/n}^{i/n}
\abs{f(x)-f(i/n)}^2\,dx\\
& \le R^2 \sum_{i=1}^n \int_{(i-1)/n}^{i/n}
\abs{x-i/n}^{2\alpha}\,dx\\
& =R^2(2\alpha+1)^{-1}n^{-2\alpha}.
\end{align*}
Consequently, asymptotic equivalence between $\E_n^1$ and $\G_n^1$
holds for any H\"older class ${\cal F}_H(\alpha,R)$ with
$\alpha>1/2$ and $R>0$ arbitrary. The approximation property of the
Haar wavelet yields even asymptotic regularity for $L^2$-Sobolev
classes of regularity $\alpha>1/2$.

For nonuniform design $0\le x_1<\cdots< x_n\le 1$ consider the same
setting as before, in particular $D_ng=(g(i/n))_i\not=(g(x_i))_i$.
We obtain for $f\in {\cal F}_H(\alpha,R)$:
\begin{align*}
&\norm{f-{\cal I}_nf}_{L^2}^2 = \sum_{i=1}^n \int_{(i-1)/n}^{i/n}
\abs{f(x)-f(x_i)}^2\,dx \le R^2 \sum_{i=1}^n \int_{(i-1)/n}^{i/n}
\abs{x-x_i}^{2\alpha}\,dx\\
&\le R^2 n^{-1}\sum_{i=1}^n \Big(n^{-1}+\abs{x_i-i/n}\Big)^{2\alpha}
\le 2R^2n^{-2\alpha}+2R^2n^{-1}\sum_{i=1}^n \abs{x_i-i/n}^{2\alpha}.
\end{align*}
By Theorem \ref{ThmMain} we have obtained the following result.

\begin{theorem}\label{ThmLocConst}
On the H\"older class ${\cal F}_H(\alpha,R)$ the Le Cam distance
between nonparametric regression with design $0<x_1^{(n)}< \cdots <
x_n^{(n)}\le 1$ and the white noise experiment satisfies
\[\Delta_{{\cal F}_H(\alpha,R)}(\E^1_n,\G^1_n)\lesssim
\sigma^{-1}R\Big(n^{1-2\alpha}+\sum_{i=1}^n
\abs{x_i^{(n)}-i/n}^{2\alpha}\Big)^{1/2}.
\]
Consequently, asymptotic equivalence holds whenever $\alpha\in
(1/2,1]$ and the design satisfies $\lim_{n\to\infty}\sum_{j=1}^n
\abs{x_i^{(n)}-i/n}^{2\alpha}=0$, e.g. if
$\max_i\abs{x_i^{(n)}-i/n}=o(n^{-1/(2\alpha)})$.
\end{theorem}

\begin{remark}
This approach does not permit to establish global equivalence for
the random design case in Section \ref{SecRanDes} because the
standard deviations of the order statistics $X_{(j)}$ decrease only
with rate $n^{-1/2}$. Treating the random design like being
equidistant yields nevertheless for estimation purposes nearly
optimal asymptotic $L^2$-risk when $\alpha>1/2$ \cite{BrownCai}.
\end{remark}

\subsection{Fourier series approximation}\label{SecFourier}

In the case of ${\cal D}=[0,1]^d$, $d\ge 1$,  and of an equidistant
design $(k/m)_{k\in \{1,\ldots,m\}^d}$ with $m=n^{1/d}\in\N$ and
odd, the Fourier system ($\iota:=\sqrt{-1}$)
\[ \phi_{ \ell}(x):=\exp(2\pi \iota\scapro{x}{\ell}),\qquad
\ell=(\ell_1,\ldots,\ell_d),\; \abs{\ell}_\infty\le \tfrac{m-1}{2},
\]
is not only $L^2$-orthonormal, but also with respect to
$\scapro{\cdot}{\cdot}_{n}$ for $D_ng:=(g(k/m))_k$:
\begin{align}
\scapro{\phi_\ell}{\phi_{\ell'}}_n&= \frac1n \sum_{k\in\{1,\ldots,m\}^d}
\phi_{\ell}(k/m)\overline{\phi_{\ell'}(k/m)}\nonumber\\
&=m^{-d}\sum_{k_1,\ldots,k_d=1}^m
\prod_{i=1}^d \exp\Big(2\pi \iota k_i(\ell_i-\ell'_i)/m\Big)\nonumber\\
&=\prod_{i=1}^d \Big(\frac 1m\sum_{\kappa=1}^m\exp\Big(2\pi \iota \kappa(\ell_i-\ell'_i)/m\Big)\Big)\nonumber\\
&=\begin{cases} 1,& \text{if } m|(\ell_i-\ell'_i)\text{ for all }i,\\
0,&\text{otherwise.}\end{cases}\label{EqFourierScaPro}
\end{align}
Consequently, the space of trigonometric polynomials
$S_n:=\spann(\phi_\ell,\,\abs{\ell}_\infty\le \tfrac{m-1}{2})$
satisfies the isometric property \eqref{EqIsometric}.

The $d$-dimensional periodic Sobolev class of regularity $s$ and
radius $R$ on $[0,1]^d$ is given by
\[ {\cal F}_{S,per}^d(s,R):=\Big\{f\in L^2([0,1]^d)\,|\,
\sum_{\ell\in\Z^d}\abs{\ell}_\infty^{2s}\abs{\scapro{f}{\phi_\ell}}^2\le
R^2\Big\}.
\]
Due to the strong cancellation property \eqref{EqFourierScaPro} of
the scalar product $\scapro{\cdot}{\cdot}_n$ we derive explicitly
\[ ({\cal I }_nf)(x)=\sum_{\abs{\ell}_\infty\le
(m-1)/2}\Big(\sum_{k\in\Z^d}\scapro{f}{\phi_{\ell+km}}\Big)\phi_\ell(x).
\]
In view of Remark \ref{Pn} we first bound the classical bias:
\[ \sup_{f\in {\cal F}_{S,per}^d(s,R)}\norm{f-P_nf}_{L^2}^2=
\sup_{f\in {\cal F}_{S,per}^d(s,R)}\sum_{\abs{\ell}_\infty\ge
(m+1)/2}\abs{\scapro{f}{\phi_\ell}}^2=R^2\Big(\frac{m+1}{2}\Big)^{-2s}.
\]
For $s>d/2$ we obtain, using the Cauchy-Schwarz inequality,
\begin{align*}
&\sup_{f\in {\cal F}_{S,per}^d(s,R)}\norm{P_nf-{\cal I
}_nf}_{L^2}^2\\
&=\sup_{f\in {\cal F}_{S,per}^d(s,R)}
\sum_{\abs{\ell}_\infty\le
(m-1)/2}\Big(\sum_{k\in\Z^d\setminus\{0\}}\scapro{f}{\phi_{\ell+km}}\Big)^2\\
&\le \Big(\sup_{f\in {\cal F}^d(s,R)}\sum_{\abs{\ell}_\infty\le
(m-1)/2}\sum_{k\in\Z^d\setminus\{0\}}\abs{\ell+km}_\infty^{2s}\scapro{f}{\phi_{\ell+km}}^2\Big)\times\\
&\quad \Big(\sup_{\abs{\ell}_\infty\le
(m-1)/2}\sum_{k\in\Z^d\setminus\{0\}}\abs{\ell+km}_\infty^{-2s}\Big)\\
&=R^2m^{-2s} \sup_{\abs{\ell}_\infty\le
(m-1)/2}\sum_{k\in\Z^d\setminus\{0\}}\abs{k+\ell/m}_\infty^{-2s}\\
&\le
R^2m^{-2s}\Big(2^{2s}(2^d-1)+\sum_{k\in\Z^d\setminus\{0\}}\abs{k}_\infty^{-2s}\Big).\\
\end{align*}
Hence, using Theorem \ref{ThmMain} we have proved
the following result, which extends the scalar results by
\cit{BrownLow} and more specifically \cit{Rohde} to any dimension
$d\ge 1$.

\begin{theorem}\label{ThmMultidim}
For $d$-dimensional periodic Sobolev classes ${\cal
F}^d_{S,per}(s,R)$ with regularity $s>d/2$ and equidistant design on
the cube $[0,1]^d$ the nonparametric regression experiment $\E_n^d$
and the Gaussian shift experiment $\G_n^d$ are asymptotically
equivalent as $n\to\infty$. The Le Cam distance satisfies
\[\Delta_{{\cal F}^d_{S,per}(s,R)}(\E^d_n,\G^d_n)\lesssim \sigma^{-1}Rn^{1/2-s/d}.\]
\end{theorem}

\section{Isomorphic approximation}\label{SecIsomorph}

\subsection{General theory}

We extend the preceding framework by merely requiring an isomorphic
property. Since it will suffice for the subsequent applications, we
specialize here immediately to $D_ng=(g(x_1),\ldots,g(x_n))$. Let
$S_n\subset {\cal L}^2({\cal D})$, $\dim S_n=n$, have the property
\begin{equation}\label{haar}
\forall g_n\in S_n:\, g_n(x_1)=\cdots=g_n(x_n)=0 \Longrightarrow
g_n=0.
\end{equation}
Let
\[ \scapro{f}{g}_n
:=\frac1n\sum_{i=1}^n f(x_i)\overline{g(x_i)}, \text{ resp.
}\scapro{v}{g}_n :=\frac1n\sum_{i=1}^n v_i\overline{g(x_i)},\quad
f,g\in {\cal L}^2,\,v\in\R^n,
\]
and $\norm{g}_n^2=\scapro{g}{g}_n$. In this notation Equation
\eqref{haar} is equivalent with the isomorphy of the norms
$\norm{\cdot}_n$ and $\norm{\cdot}_{L^2}$ on $S_n$:
\begin{equation}\label{EqIsomorph}
\exists\,
A_n,B_n>0\;\forall\,g_n\in S_n:\; A_n\norm{g_n}_{L^2}\le
\norm{g_n}_n\le B_n\norm{g_n}_{L^2}.
\end{equation}
We choose any $L^2$-orthonormal basis $(\phi_j)_{1\le j\le n}$ of
$S_n$ and introduce the linear mappings
\[ \Pi_n,\,{\cal I }_n:{\cal
L}^2({\cal D}) \to S_n,\quad
\Pi_ng:=\sum_{j=1}^n\scapro{g}{\phi_j}_n\phi_j,\quad {\cal
I}_ng:=\Pi_n|_{S_n}^{-1}\Pi_ng.
\]
Observe the following properties: for $g_n,h_n\in S_n$ we have
$\scapro{\Pi_ng_n}{h_n}=\scapro{g_n}{h_n}_n$ and
$\norm{\Pi_n|_{S_n}}\le B_n$, $\norm{(\Pi_n|_{S_n})^{-1}}\le
A_n^{-1}$; ${\cal I }_n$ is a projection onto $S_n$ and ${\cal
I}_ng$ interpolates $g$ at the design points $(x_i)$; $\Pi_n$ and
${\cal I}_n$ are independent of the choice of the basis $(\phi_j)$.

The regression experiment $\E^d_n$ can be transformed to a
functional Gaussian shift by expanding the observations $(Y_i)$ in
the basis $(\phi_j)$:
\begin{equation}\label{EqZ1}
Z_1:=\sum_{j=1}^n\scapro{Y}{\phi_j}_n\phi_j=\Pi_nf+\frac{\sigma}{\sqrt{n}}(\Pi_n|_{S_n})^{1/2}\zeta\quad\in
S_n,
\end{equation}
with Gaussian white noise
$\zeta:=(\Pi_n|_{S_n})^{-1/2}(\sqrt{n}\sum_{j=1}^n\scapro{\eps_j}{\phi_j}_n\phi_j)\sim
{\cal N}(0,\Id_{S_n})$ because
\[
E[\scapro{\zeta}{g_n}\scapro{\zeta}{h_n}] =
\scapro{(\Pi_n|_{S_n})^{-1/2}g_n}{(\Pi_n|_{S_n})^{-1/2}h_n}_{n}=\scapro{g_n}{h_n},\quad
g_n,h_n\in S_n.
\]
By applying $(\Pi_n|_{S_n})^{-1/2}$ and $(\Pi_n|_{S_n})^{-1}$,
respectively, we conclude that the regression experiment $\E^d_n$ is
also equivalent to observing
\begin{align}
Z_2&=(\Pi_n|_{S_n})^{-1/2} Z_1=(\Pi_n|_{S_n})^{1/2}{\cal I
}_nf+\frac{\sigma}{\sqrt{n}}\zeta\quad \in S_n,\label{EqZ2}\\
Z_3&=(\Pi_n|_{S_n})^{-1} Z_1={\cal I
}_nf+\frac{\sigma}{\sqrt{n}}(\Pi_n|_{S_n})^{-1/2}\zeta\quad \in
S_n\label{EqZ3}
\end{align}
with $\zeta\sim{\cal N}(0,\Id_{S_n})$.

\begin{theorem}\label{ThmMainlocal}
The regression experiment $\E^d$ is equivalent with each of the
experiments given by observing $Z_1$ in \eqref{EqZ1}, $Z_2$ in
\eqref{EqZ2} and $Z_3$ in \eqref{EqZ3}, respectively.

The Le Cam distance between $\E^d_n$ and $\G^d_n$ for the class
${\cal F}^d$ satisfies the bounds
\begin{align}
\Delta_{{\cal F}^d}(\E^d_n,\G^d_n)&\le
1-2\Phi\Big(-\frac{\sqrt{n}}{2\sigma}\sup_{f\in{\cal F}^d}
\norm{f-\Pi_n|_{S_n}^{1/2}{\cal I
}_nf}_{L^2}\Big),\label{EqDeltaZ2}\\
\Delta_{{\cal F}^d}(\E^d_n,\G^d_n)&\le
1-2\Phi\Big(-\frac{\sqrt{n}}{2\sigma}\sup_{f\in{\cal F}^d}
\norm{f-{\cal I
}_nf}_{L^2}\Big)+\sqrt2\norm{(\Pi_n|_{S_n})^{-1}-\Id_{S_n}}_{HS},\label{EqDeltaZ3}
\end{align}
where $\norm{\cdot}_{HS}$ denotes the Hilbert-Schmidt norm of an
operator.
\end{theorem}

\begin{proof}
It remains to prove the second part. The first bound
\eqref{EqDeltaZ2} follows from the equivalence with observing $Z_2$
by the same arguments as for Theorem \ref{ThmMain}. To establish
\eqref{EqDeltaZ3}, we use the fact that the Hellinger distance
between two multivariate normal distributions with the same mean
satisfies
\begin{equation}\label{EqHellBound}
H^2(N(\mu,\alpha\Sigma),N(\mu,\alpha\Id_{\R^n}))\le
2\norm{\Sigma-\Id_{\R^n}}_{HS}^2,\quad \Sigma\in\R^{n\times
n},\,\alpha>0,
\end{equation}
which follows e.g. from \cite[Lemma 3]{Brownetal} via the
diagonalisation $\Sigma=O^\top \diag(\lambda_1,\ldots,\lambda_n)O$
and the property
$\norm{\Sigma-\Id_{\R^n}}_{HS}^2=\norm{O(\Sigma-\Id_{\R^n})O^\top}_{HS}^2=\sum_{i=1}^n\lambda_i^2$.
Therefore the total variation distance between the laws of $Z_3$ and
$Z_4:={\cal I }_nf+\frac{\sigma}{\sqrt{n}}\zeta$ is bounded by
\[ \norm{{\cal L}(Z_3)-{\cal L}(Z_4)}_{TV}\le H({\cal L}(Z_3),{\cal
L}(Z_4))\le \sqrt2\norm{(\Pi_n|_{S_n})^{-1}-\Id_{S_n}}_{HS}.
\]
The by now standard arguments yield with obvious notation
\begin{align*}
\Delta_{{\cal F}^d}(\E_n^d,\G_n^d)&=\Delta_{{\cal
F}^d}(Z_3,\G_n^d)\le \Delta_{{\cal F}^d}(Z_4,\G_n^d)+\Delta_{{\cal
F}^d}(Z_4,Z_3)\\
&\le 1-2\Phi\Big(-\frac{\sqrt{n}}{2\sigma}\sup_{f\in{\cal F}^d}
\norm{f-{\cal I
}_nf}_{L^2}\Big)+\sqrt2\norm{(\Pi_n|_{S_n})^{-1}-\Id_{S_n}}_{HS},
\end{align*}
as asserted.
\end{proof}

\subsection{Linear spline approximation}

Let us briefly expose how the approach by \cit{Carter} fits into the
isomorphic framework. As in Section \ref{SecFourier}, we consider
equidistant design points $(k/m)_{k\in \{1,\ldots,m\}^d}$ with
$m=n^{1/d}\in\N$ and periodic functions on the unit cube ${\cal
D}=[0,1]^d$. The space $S_n$ is spanned by the periodized and
tensorized linear B-splines
\[ b_k(x)=b_k(x_1,\ldots,x_d)=\prod_{r=1}^d \bar b(mx_r-k_r \mod
1),\quad \bar b:={\bf 1}_{[-1/2,1/2]}\ast {\bf 1}_{[-1/2,1/2]},
\]
indexed by $k\in \{1,\ldots,m\}^d$. For $\alpha\in (1,2]$ it is well
known (cf. \cit{DeBoor}) that interpolation on $S_n$ for the
periodic H\"older class
\[ {\cal F}^d_{H,per}(\alpha,R):=\big\{ f\in C^\alpha(\R^d)\,\big|\, f\text{ $\Z^d$-periodic, }\sup_{x\not=
y}\abs{\nabla f(x)-\nabla f(y)}/\abs{x-y}^{\alpha-1}\le R\big\}
\]
satisfies the estimate
\begin{equation}\label{InEst}
\sup_{f\in{\cal F}_{H,per}^d(\alpha,R)} \norm{f-{\cal
I}_nf}_{L^2([0,1]^d)}\lesssim Rn^{-\alpha/d}.
\end{equation}
On the other hand, we have for $g_n\in S_n$
\[ \norm{g_n}_{L^2}^2=\Big\lVert\sum_{k\in
\{1,\ldots,m\}^d} g_n(k/m)b_k\Big\rVert_{L^2}^2=\sum_{k,\ell\in
\{1,\ldots,m\}^d} \scapro{b_k}{b_\ell}_{L^2}g_n(k/m)g_n(\ell/m)
\]
with $\scapro{b_k}{b_\ell}_{L^2}=0$ for $\abs{k-\ell}_\infty>1$ and
$\scapro{b_k}{b_\ell}_{L^2}=4^{\#\{r:k_r=\ell_r\}}/(6^d n)$ for
$\abs{k-\ell}_\infty\le 1$. Since
$\sum_{\ell}\scapro{b_k}{b_\ell}_{L^2}=\scapro{b_k}{1}_{L^2}=n^{-1}$,
a weighted Cauchy-Schwarz inequality yields
\[ \norm{g_n}_{L^2}^2\le n^{-1}\sum_{k\in
\{1,\ldots,m\}^d}
g_n(k/m)^2=\scapro{g_n}{g_n}_{n}=\scapro{\Pi_ng_n}{g_n}_{L^2}
\]
and we conclude, using the ordering of symmetric operators, that
$(\Pi_n|_{S_n})^{-1}\le \Id_{S_n}$.

Adding to the observation $Z_3$ in \eqref{EqZ3} independent Gaussian
noise $\eta\sim {\cal
N}(0,\frac{\sigma^2}{n}(\Id_{S_n}-(\Pi_n|_{S_n})^{-1}))$, we infer
that the regression experiment $\E_n^d$ is more informative than
observing
\begin{equation}\label{EqZ5}
Z_5:=Z_3+\eta={\cal I}_nf+\frac{\sigma}{\sqrt{n}}\tilde\zeta\quad\in
S_n
\end{equation}
with Gaussian white noise
$\tilde\zeta:=(\Pi_n|_{S_n})^{-1/2}\zeta+n^{1/2}\sigma^{-1}\eta\sim{\cal
N}(0,\Id_{S_n})$. This randomization together with estimate
\eqref{InEst} shows that the regression experiment $\E_n^d$ is
asymptotically at least as informative as the Gaussian experiment
$\G_n^d$ on H\"older classes ${\cal F}_{per}^d(\alpha,R)$ with
$\alpha>d/2$ and $d\in\{1,2,3\}$. Together with an (easier)
randomization in the other direction and a more sophisticated
boundary treatment for non-periodic function classes this reproduces
the proof by \cit{Carter} for asymptotic equivalence of regression
and white noise experiments in dimensions 2 and 3. For B-splines of
higher order the interpolation property $b_k(i/m)=\delta_{k,i}$ gets
lost and $(\Pi_n|_{S_n})^{-1}\le \Id_{S_n}$ cannot be shown such
that a more refined analysis is needed. This will be accomplished in
the next section for a similar approach using compactly supported
wavelets.

\subsection{Wavelet multiresolution analysis}

{\bf The construction.} Let us assume an equidistant dyadic design
$(k2^{-j})_{k\in \{1,\ldots,2^j\}^d}$ with $n=2^{dj}$ points for
some $j\in\N$ and ${\cal D}=[0,1]^d$. We consider a wavelet
multiresolution analysis $(V_j)_{j\ge 0}$ on $L^2([0,1]^d)$ obtained
by periodisation and tensor products. Let $\bar\phi$ be a standard
orthonormal scaling function of an $r$-regular multiresolution
analysis for $L^2(\R)$, that is $(\bar\phi(\cdot+k))_{k\in\Z}$ forms
an orthonormal system in $L^2(\R)$ and satisfies $\int\bar\phi=1$ as
well as the polynomial exactness condition that
$\sum_{k\in\Z}k^q\bar\phi(x-k)-x^q$ is a polynomial of maximal
degree $q-1$ for all $q=0,\ldots,R-1$ \cite[Thm. 16.1]{Cohen}. We
suppose that $\bar\phi$ has compact support in $[-S+1,S]$, like in
Daubechies's construction, so that the functions
$\phi_{jk}:[0,1]^d\to\R$, $j\ge 1$, $k\in \{1,\ldots,2^j\}^d$, with
\[
\phi_{jk}(x_1,\ldots,x_d):=\sum_{m\in\Z^d}2^{jd/2}\prod_{i=1}^d\bar\phi(2^jx_i-k_i+2^jm_i)
\]
are well defined and form an orthonormal system in $L^2([0,1]^d)$
\cite[Prop. 2.21]{Wojtaszczyk}. We set
$S_{2^{jd}}:=V_j:=\spann\{\phi_{jk}\,|\,k\in\{1,\ldots,2^j\}^d\}$.\\

{\bf\noindent Periodic approximation.} Polynomial exactness and
continuity of $\bar\phi$ imply for $q=0,\ldots,R-1$ and any $x\in\R$
\cite{SwePie}
\[ \sum_{m\in\Z}(x+m)^q\bar\phi(x+m)=\int_{-\infty}^\infty
x^q\bar\phi(x)\,dx.
\]
This identity is fundamental for our purposes because it implies for
$\Z^d$-periodic functions $h:\R^d\to\R$ that coincide with a
polynomial $p$ of maximal degree $R-1$ on
$\prod_{i=1}^d[2^{-j}(k_i-S-1),2^{-j}(k_i+S)]$:
\begin{align*}
\scapro{h}{\phi_{jk}}_{L^2}&=\sum_{m\in\Z^d}2^{jd/2}\int_{[0,1]^d}h(x)\prod_{i=1}^d\bar\phi(2^j(x_i+m_i)-k_i)\,dx\\
&=2^{jd/2}\int_{\R^d}h(x)\prod_{i=1}^d\bar\phi(2^jx_i-k_i)\,dx\\
&=2^{-jd/2}\int_{[-S-1,S]^d}p(2^{-j}(x+k))\prod_{i=1}^d\bar\phi(x_i)\,dx\\
&=2^{-jd/2}\sum_{m\in\Z^d}p(2^{-j}(m+k))\prod_{i=1}^d\bar\phi(m_i)\\
&=2^{-jd/2}\sum_{m\in\{1,\ldots,2^j\}^d}h(2^{-j}m)\phi_{jk}(2^{-j}m)\\
&=n^{1/2}\scapro{h}{\phi_{jk}}_n,
\end{align*}
where we identified $n=2^{jd}$. For any $\Z^d$-periodic function
$g\in H^s_{S,per}([0,1]^d)$ with $s\in (d/2,R)$ this local
polynomial reproduction property implies by standard, but
sophisticated arguments for direct estimates \cite[Thm. 30.6]{Cohen}
\begin{equation}\label{EqApproxCohen}
\norm{g-\Pi_ng}_{L^2} \lesssim
2^{-js}\norm{g}_{H^s}=n^{-s/d}\norm{g}_{H^s},
\end{equation}
where $\norm{\cdot}_{H^s}$ denotes the standard $L^2$-Sobolev norm
of regularity $s$ on $[0,1]^d$. We split the bias term and obtain by
functional calculus
\begin{align*}
\norm{f-\Pi_n|_{S_n}^{-1/2}\Pi_nf}_{L^2} &\le
\norm{f-\Pi_nf}_{L^2}+\norm{\Pi_nf-\Pi_n|_{S_n}^{-1/2}\Pi_nf}_{L^2}\\
&=
\norm{f-\Pi_nf}_{L^2}+\norm{h(\Pi_n|_{S_n})(\Id-\Pi_n)\Pi_nf}_{L^2}
\end{align*}
with $h:\R^+\to\R$, $h(x):=1/(x+x^{1/2})=(x^{-1/2}-1)/(1-x)$. Since
$h$ decreases monotonically and $h(x)\le x^{-1/2}$, we have
$\norm{h(\Pi_n|_{S_n})}_{L^2\to L^2}\le \lambda_{min}^{-1/2}$ with
the smallest eigenvalue $\lambda_{min}$ of $\Pi_n|_{S_n}$.

$\Pi_n|_{S_n}$ satisfies for $n=2^{jd}\ge 2S-1$ the following
scaling property:
\begin{align*}
\scapro{\Pi_n\phi_{jk}}{\phi_{j\ell}}_{L^2}
&=\frac1n\sum_{\nu\in\{1,\ldots,2^j\}^d}\phi_{jk}(\nu2^{-j})\phi_{j\ell}(\nu2^{-j})\\
&=\sum_{m\in\Z^d}\sum_{\nu\in\{1,\ldots,2^j\}^d}
\prod_{a=1}^d\Big(\bar\phi((\nu-k+2^jm)_a)\bar\phi((\nu-\ell+2^jm)_a)\Big)\\
&=\prod_{a=1}^d\Big(\sum_{b\in\Z}
\bar\phi(b-k_a)\bar\phi(b-\ell_a)\Big).
\end{align*}
Since $\bar\phi$ has compact support, the series is just a finite
sum and $\Pi_n$ has a bounded Toeplitz matrix representation in
terms of $(\phi_{jk})$. Using Fourier multipliers it follows that
$\scapro{\Pi_ng_n}{g_n}_{L^2}\ge A_{\bar\phi}^2\norm{g_n}_{L^2}^2$,
$g_n\in S_n$, with $A_{\bar\phi}:=\inf_{u\in[0,2\pi]}
\abs{\sum_{k\in\Z}\bar\phi(k)e^{\iota ku}}^{d}$, independently of
$n$. Due to the compact support of $\bar\phi$, we have
$A_{\bar\phi}>0$ iff the trigonometric polynomial
$\sum_{k\in\Z}\bar\phi(k)e^{\iota ku}$, $u\in [0,2\pi]$, does not
vanish. It is well known \cite[Lemma 3]{SwePie} that this is exactly
the condition to ensure that the multiresolution analysis is also
generated by an interpolating scaling function. It can be checked
for standard Daubechies scaling functions, e.g. by showing
$\abs{\bar\phi(k_0)}>\sum_{k'\not=k_0}\abs{\bar\phi(k')}$ for some
$k_0\in\Z$. Moreover, gaining more flexibility by considering the
shifted spaces based on $\bar\phi_\tau=\bar\phi(\cdot-\tau)$,
$\tau\in (0,1)$, a wavelet multiresolution analysis will almost
always satisfy $A_{\bar\phi_\tau}>0$ for some value of $\tau$, cf.
\cit{SwePie} and the references therein.

We arrive at
\[ \norm{f-\Pi_n|_{S_n}^{-1/2}\Pi_nf}_{L^2}\le
\norm{f-\Pi_nf}_{L^2}+A_{\bar\phi}^{-1/2}\norm{(\Id-\Pi_n)\Pi_nf}_{L^2}.
\]
Because of $\norm{\Pi_nf}_{H^s}\to\norm{f}_{H^s}$ \cite[Thm.
30.7]{Cohen} we derive from \eqref{EqApproxCohen} the uniform
estimate over $f\in {\cal F}_{S,per}^d(s,R)$
\[ \norm{f-\Pi_n|_{S_n}^{-1/2}\Pi_nf}_{L^2}\le \norm{f-\Pi_nf}_{L^2}+A_{\bar\phi}^{-1/2}
\norm{(\Id-\Pi_n)\Pi_nf}_{L^2}\lesssim Rn^{-s/d}.
\]
Hence, the estimate in \eqref{EqDeltaZ2} yields asymptotic
equivalence between the regression and the white noise experiment
for any class ${\cal F}_{S,per}^d(s,R)$ with $s>d/2$.

This result provides another way for constructing explicitly the
transformation between the regression and the white noise setting.
It has no more theoretical implications than the Fourier basis
approach, but it paves the way for proving asymptotic equivalence
for non-periodic function classes.\\

{\bf\noindent Non-periodic approximation.} Since every $\phi_{jk}$
has support length $2^{-j}(2S-1)$, only those functions $\phi_{jk}$
with $k_r\in\{1,\ldots,S-2\}\cup\{2^j-S+1,\ldots,2^j\}$ for some
$r=1,\ldots,d$ cross the boundary and are periodized at all.
Therefore, the same derivation using only interior scaling functions
shows that the regression experiment $\E^d_n$ for the general
Sobolev function class
\[ {\cal F}_S^d(s,R):=\{f\in H^s([0,1]^d)\,|\,\norm{f}_{H^s}\le R\}\]
is asymptotically more informative than the restricted white noise
experiment $\bar\G^d_{n}$ given by observing
\begin{equation}\label{EqGrestr}
dY(x)=f(x)\,dx+\frac{\sigma}{\sqrt{n}}dB(x),\quad x\in
[\delta_n,1-\delta_n]^d\text{ with }\delta_n:=(2S-1)n^{-1/d}.
\end{equation}
Although $\bar\G^d_{n}$ is a priori less informative than $\G^d_n$,
we may use classical extrapolation, e.g. the Taylor polynomial
$T_f^y$ of order $\floor{s}$ around $y\in [\delta_n,1-\delta_n]^d$.
We define at the points $x\in [0,1]^d\setminus
[\delta_n,1-\delta_n]^d$ the extrapolation $\tilde
f(x)=T_f^{y_x}(x)$ for a point $y_x\in [\delta_n,1-\delta_n]^d$ with
$\abs{y_x-x}_\infty\le 2\delta_n$, selected in a measurable way, and
$\tilde{f}(x)=f(x)$ otherwise. We thereby achieve
\[ \Big(\int_{[0,1]^d}\abs{\tilde f(x)-f(x)}^2\,dx\Big)^{1/2}\lesssim
R n^{-s/d}
\]
such that
\[ \Delta_{{\cal F}_S^d(s,R)}(\bar\G^d_{n},\G^d_n)\lesssim
\sigma^{-1}R n^{1/2-s/d}.
\]
This means that $\bar\G^d_{n}$ and $\G^d_n$ are asymptotically
equivalent for $s>d/2$ and we have obtained a result for function
classes without periodicity condition.

\begin{theorem}\label{ThmMultidimWav}
For general $d$-dimensional Sobolev classes ${\cal F}^d_S(s,R)$ with
regularity $s>d/2$ and equidistant design on the cube $[0,1]^d$ the
nonparametric regression experiment $\E_n^d$ and the Gaussian white
noise experiment $\G_n^d$ are asymptotically equivalent as
$n\to\infty$. The Le Cam distance satisfies
\[\Delta_{{\cal F}^d_S(s,R)}(\E^d_n,\G^d_n)\lesssim \sigma^{-1}Rn^{1/2-s/d}.
\]
\end{theorem}

\noindent{\bf Discussion.} The property that a wavelet estimator
based on an equidistant regression model and a corresponding
estimator based on a white noise model are asymptotically close is
well known, see e.g. \cit{DonJoh2} and \cit{JohnSilv}.
Interestingly, both papers show identical asymptotics of the
$L^2$-risk for standard estimators uniformly over balls in Besov
spaces $B^s_{p,q}([0,1])$ with $s>1/p$ or $s=p=1$. Since $B^s_{p,q}$
embeds into the Sobolev space $H^\sigma$ for $s>\sigma$ and
$s-1/p>\sigma-1/2$, Theorem \ref{ThmMultidimWav} provides more
generally asymptotic equivalence for Besov classes with $s>1/p$ and
$p<2$. The counterexample in \cit{BrownLow} shows, however, that for
$s\le 1/2$ and all $p\in [1,\infty]$, asymptotic equivalence breaks
down. Similarly, if $\psi\in B^1_{1,1}$ is a function with support
in $(0,1)$ and $\norm{\psi}_{L^2}=1$, then $\psi_{n}(x):=\psi(nx)$
has support in $(0,1/n)$, $L^2$-norm $\norm{\psi_n}_{L^2}=n^{-1/2}$
and Besov norm $\norm{\psi_n}_{B^1_{1,1}}\thicksim 1$. Hence,
testing the signal $f=0$ versus $f=\psi_n$ has nontrivial power in
the white noise model $\G_n^1$, while both signals generate exactly
the same observations in the regression model $\E_n^1$. We conclude
that $\G_n^1$ and $\E_n^1$ are not asymptotically equivalent on
Besov classes with $s=1$, $p=1$. An intriguing example for the
important class of bounded variation functions is given by
$\psi_n(x)=\sqrt{2}{\bf 1}_{[1/4n,3/4n]}(x)$. Asymptotic equivalence
between Gaussian regression and white noise is indeed an
$L^2$-theory and we cannot gain by measuring smoothness in an
$L^p$-sense, $p\not=2$.

Let us also mention that the (asymptotically negligible) loss in
information due to neglecting boundary coefficients in the
construction seems unavoidable. The wavelets on an interval
\cite{CDV} use nonorthogonal boundary corrections and can therefore
not be used, while the coiflet approach by \cit{JohnSilv} also
involves some information loss at the boundary, cf. their remark on
dimensions before Proposition 2.

\section{Random design}\label{SecRanDes}

\subsection{The general idea}

Denote by $U([0,1]^d)$ the uniform distribution on the cube ${\cal
D}=[0,1]^d$.

\begin{definition}
Let $\E^d_{n,r}$ be the compound experiment obtained from observing
independent random design points $X_i\sim U([0,1]^d)$,
$i=1,\ldots,n$, and the regression
\[ Y_i=f(X_i)+\sigma\eps_i,\qquad i=1,\ldots,n,
\]
for $n\in\N$ and $f:[0,1]^d\to\R$ in some class ${\cal F}^d\subset
{\cal L}^2([0,1]^d)$ and with i.i.d. random variables $\eps_i\sim
{\cal N}(0,1)$, independent of the design.
\end{definition}

We place ourselves into the isomorphic setting, that is we are given
an $L^2([0,1]^d)$-orthonormal basis $(\phi_j)_{j\ge 1}$ and we set
$S_n=\spann(\phi_1,\ldots,\phi_n)$. For the moment we merely assume
that $S_n$ is chosen to satisfy the isomorphic condition
\eqref{haar} given the random design points $(X_i)_{1\le i\le n}$.
Later, certain parts will rely on fine properties of the Fourier
basis. Conditionally on the design the regression experiment is
equivalent to observing
\[ Z_1:=\sum_{j=1}^{n}\scapro{Y}{\phi_j}_n\phi_j
=\Pi_{n}f+\frac{\sigma}{\sqrt{n}}(\Pi_n|_{S_n})^{1/2}\zeta\quad\in
S_n
\]
with white noise $\zeta\sim N(0,\Id_{S_n})$. Let us briefly comment
why the foregoing approaches using $Z_2$ in \eqref{EqZ2} or $Z_3$ in
\eqref{EqZ3} will not succeed here. For
$Z_2=(\Pi_n|_{S_n})^{-1/2}Z_1$ we need to have
$\norm{(\Pi_n|_{S_n}^{1/2}-\Id){\cal I}_nf}_{L^2}$ and $\norm{{\cal
I}_nf-f}_{L^2}$ of smaller order than $n^{-1/2}$. The second
property can be ensured for Sobolev classes of regularity $s>d/2$ as
before. The first property, however, will not hold. By empirical
process theory, we have for $g_1,g_2\in S_n$ approximately
$\scapro{\Pi_{n}g_1}{g_2}_{L^2}=\scapro{g_1}{g_2}_n\approx
\scapro{g_1}{g_2}_{L^2}+n^{-1/2}\int g_1g_2dB^0$ with a Brownian
bridge $B^0$. By the linearisation $(1+h)^{1/2}-1\approx h/2$ and
taking expectation with respect to the random design, we find
\begin{align*}
E\Big[\norm{((\Pi_n|_{S_n})^{1/2}-\Id){\cal
I}_nf}_{L^2}^2\Big]&\thicksim E\Big[ \sum_{j=1}^n \babs{n^{-1/2}\int
({\cal I}_n f) \phi_j dB^0}^2\Big]\\
&\thicksim n^{-1}\sum_{j=1}^n \int \abs{\phi_j}^2\abs{{\cal I}_n
f}^2 \gtrsim \norm{{\cal I}_nf}_{L^2}^2.
\end{align*}
Hence, in the mean over the random design this term does not tend to
zero. When considering $Z_3=(\Pi_n|_{S_n})^{-1}Z_1$, we would need
$\norm{(\Pi_n|_{S_n})^{-1}-\Id_{S_n}}_{HS}\to 0$, compare Bound
\eqref{EqDeltaZ3}, but the mean over this term is by the same
approximations of order $n$. The main defect in these approaches is
that we do not take advantage of the regularity of $f$.

The new idea is based on a two-level procedure, which can be
interpreted as a localisation approach, cf. \cit{Nussbaum}. We
choose an intermediate level $n_0<n$ and split
$S_n=S_{n_0}+U_{n_0}^n$ with the $\norm{\cdot}_n$-orthogonal
complement $U_{n_0}^n$ of $S_{n_0}$ in $S_n$. On the low-frequency
space $S_{n_0}$ we use the empirical orthogonal projection
$P^n_{n_0}Y$ of the data onto $S_{n_0}$. This construction is
analogous to $Z_3$ in \eqref{EqZ3} and the heteroskedasticity in the
noise term will become asymptotically negligible provided
$n_0=o(n^{1/2})$.

On the high-frequency part $U_{n_0}^n$ of $S_n$ we transform to a
Gaussian shift with white noise, which is independent of the noise
in $S_{n_0}$, in the spirit of $Z_2$ in \eqref{EqZ2}. In order to
take advantage of the regularity of $f$, however, we do not use the
standard square root operator $\Pi_n^{-1/2}$ to whiten the noise,
but the adjoint $T^\ast$ of an operator $T:S_n\to S_n$ which has an
upper triangular matrix representation in the basis $(\phi_j)$ and
satisfies $TT^\ast=(\Pi_n|_{S_n})^{-1}$ (as in the Cholesky
decomposition). Since $T^\ast$ is a unitary transformation of
$(\Pi_n|_{S_n})^{-1/2}$, the noise part remains white. Due to the
triangular structure, the signal coefficients
$\scapro{T^\ast\Pi_nf}{\phi_j}_{L^2}=\scapro{T^{-1}{\cal
I}_nf}{\phi_j}_{L^2}$ do not involve the (usually large)
coefficients $\scapro{{\cal I}_nf}{\phi_k}_{L^2}$ for indices $k$
smaller than $j$. Moreover, for the Fourier basis the other
off-diagonal matrix entries of $T^{-1}$ are centred and
uncorrelated, while the deviations in the diagonal entries grow with
the frequencies, but are exactly counter-balanced by the decay of
the Fourier coefficients for Sobolev function classes. Provided
$n_0\to\infty$, this high-frequency transformation will imply
asymptotic equivalence.

\subsection{The main result}

Let us specify the transformation $T$ concretely based on the
Gram-Schmidt procedure for orthonormalisation with respect to
$\norm{\cdot}_n$. For $j\le n$ denote by $P_j,\,P_j^n:S_n\to S_n$
the $L^2$-orthogonal and $\norm{\cdot}_n$-orthogonal projections
onto $S_j$, respectively, and set $P_0^n:=0$. We obtain an
$\norm{\cdot}_n$-orthonormal basis $(\phi_j^n)$ of $S_n$ via
\[
\phi_j^n:=\frac{\phi_j-P_{j-1}^n\phi_j}{\norm{\phi_j-P_{j-1}^n\phi_j}_n},\quad
j=1,\ldots,n.
\]
Then $\phi_j^n$ is in $S_j$ and the $\norm{\cdot}_n$-orthogonality
$\phi_j^n\perp_n S_{j-1}$ holds. Defining $T:S_n\to S_n$ via
$T\phi_j:=\phi_j^n$, we see that $T$ satisfies
$\scapro{T\phi_{j'}}{\phi_{j}}_{L^2}=0$ for $j>j'$ and is an
isometry between $(S_n,\norm{\cdot}_{L^2})$ and
$(S_n,\norm{\cdot}_n)$ such that $\Pi_n|_{S_n}=(TT^\ast)^{-1}$.
The noise terms $(\scapro{\eps}{\phi_j^n}_n)_{1\le j\le n}\sim {\cal
N}(0,n^{-1})$ are therefore independent and
\[
P_{n_0}^n\eps:=\sum_{j=1}^{n_0}\scapro{\eps}{\phi_j^n}_n\phi_j^n
=\sum_{j=1}^{n_0}\scapro{\eps}{\phi_j^n}_nT\phi_j\sim {\cal
N}(0,n^{-1}T|_{S_{n_0}}T|_{S_{n_0}}^\ast).
\]
Using $T|_{S_{n_0}}T|_{S_{n_0}}^\ast=(\Pi_n|_{S_{n_0}})^{-1}$, we
introduce the rescaled covariance operator $\Sigma:S_n\to S_n$ via
\[ \Sigma
g_n:=(\Pi_n|_{S_{n_0}})^{-1}P_{n_0}g_n+(\Id_{S_n}-P_{n_0})g_n,\quad
g_n\in S_n.
\]
The regression experiment is then transformed to
observing
\begin{align}
Z_r&:=\sum_{j=1}^{n_0}
\scapro{Y}{\phi_j^n}_n\phi_j^n+\sum_{j=n_0+1}^n
\scapro{Y}{\phi_j^n}_n\phi_j\quad\in S_n\label{EqZr}\\
&=P_{n_0}^nf+T^{-1}(P^n_n-P^n_{n_0})f+n^{-1/2}\sigma\Sigma^{1/2}\zeta\quad\in
S_n\nonumber
\end{align}
with Gaussian white noise $\zeta\sim N(0,\Id_{S_n})$, conditional on
the random design.

\begin{example}
Let us consider the Haar basis. Write
$I_{jk}=[2^{-j}k,2^{-j}(k+1))$, $N_{jk}=\#\{i:X_i\in I_{jk}\}$ and
$\psi_{jk}=2^{j/2}({\bf 1}_{I_{j+1,2k}}-{\bf 1}_{I_{j+1,2k+1}})$ for
$j\ge 0$, $k=0,\ldots, 2^j-1$. By construction the transformed basis
function $\psi_{jk}^n$ has support $I_{jk}$, is constant on
$I_{j+1,2k}$, $I_{j+1,2k+1}$ and satisfies
$\scapro{\psi_{jk}^n}{{\bf 1}_{I_{jk}}}_n=0$,
$\norm{\psi_{jk}^n}_n=1$. We infer
\[
\psi_{jk}^n=C_{jk}\big(N_{j,2k}^{-1}{\bf
1}_{I_{j,2k}}-N_{j,2k+1}^{-1}{\bf 1}_{I_{j,2k+1}}\big), \quad
C_{jk}^2=nN_{j+1,2k}N_{j+1,2k+1}/N_{jk}.
\]
This is exactly the application of our framework underlying previous
one-dimensional constructions \cite[Eq. (2.8)]{Brownetal}. Because
here $S_n$ is for most design realisations not isomorphic,
additional randomisations are needed.
\end{example}

For the following general $d$-dimensional theorem we consider the
construction \eqref{EqZr} in terms of the Fourier basis functions
$\phi_j(x)=\exp(2\pi \iota \scapro{\ell(j)}{x})$ with an enumeration
$\ell:\N\to\Z^d$ of $\Z^d$ satisfying $\abs{\ell(j)}_{\ell^2}\le
\abs{\ell(j')}_{\ell^2}$ for $j\le j'$ (i.e. sorted in the order of
magnitudes of the frequencies).

\begin{theorem}\label{ThmRanDes}
For $d$-dimensional periodic Sobolev classes ${\cal
F}^d_{S,per}(s,R)$ with regularity $s>d/2$ the nonparametric
regression experiment $\E_{n,r}^d$ with random design and the
Gaussian shift experiment $\G_n^d$ are asymptotically equivalent as
$n_0,n\to\infty$ and $n_0=o(n^{1/2})$. The Le Cam distance satisfies
\[\Delta_{{\cal F}^d_{S,per}(s,R)}(\E^d_{n,r},\G^d_n)\lesssim n^{-1/2}n_0+\sigma^{-1}Rn_0^{1/2-s/d}.
\]
\end{theorem}

\begin{remark}
The asymptotically optimal choice of $n_0$ is given by $n_0\thicksim
n^{d/(2s+d)}$, which yields a bound on the Le Cam distance of order
$n^{(d-2s)/(2d+4s)}$. Note that this choice $n_0\thicksim
n^{d/(2s+d)}$ corresponds exactly to the optimal dimension of the
approximation spaces in nonparametric regression and is also used by
\cit{Gaiffas} for his two-level construction of optimal confidence
bands.
\end{remark}

\begin{proof}
In order to bound the Le Cam distance for compound experiments, we
use that for distributions $K\otimes P$ and $K'\otimes P$, defined
on $(\Omega\times\Omega',{\cal F}\otimes{\cal F}')$ by the measure
$P$ on $\cal F$ and the Markov kernels $K,K'$ from $\Omega$ to
${\cal F}'$, the total variation distance can be calculated by
conditioning:
\[ \norm{K\otimes P-K'\otimes
P}_{TV({\cal F}\otimes{\cal F'})}
=\int\norm{K(\omega,\cdot)-K'(\omega,\cdot)}_{TV({\cal
F'})}\,P(d\omega).
\]
Therefore we can first work conditionally on the design and then
take expectations for $(X_i)$. Moreover, the white noise experiment
$\G_n^d$ is equivalent to the compound experiment of $\G_n^d$ and
the observation of the random design points because the latter is a
trivial randomisation of $\G_n^d$.

It is a remarkable property of the Fourier basis that $S_n$ is
almost surely isomorphic, cf. Theorem 1.1 in \cit{BassGro}. In
Proposition \ref{PropOmegaj} below we prove that the event
\begin{equation}\label{EqOmegaj}
\Omega_j^n:=\{\forall g\in S_j:\;\tfrac12\norm{g}_{L^2}\le
\norm{g}_n\le 2\norm{g}_{L^2}\}
\end{equation}
for $j\log(j)=o(n)$ even satisfies $P((\Omega_{j}^n)^\complement)\to
0$ with a convergence rate faster than any polynomial in $n$. This
is much tighter with respect to the subspace dimension than what can
be derived from \cit{BassGro}. In order to establish asymptotic
equivalence, it suffices therefore to estimate the total variation
distances on the event $\Omega_{n_0}^n$.

By \eqref{EqZr}, the regression experiment $\E_{n,r}^d$ is
equivalent to observing $Z_r$ together with the design. Introducing
\begin{equation}\label{EqZr'}
\bar Z_r:=P_nf+\sigma n^{-1/2}\zeta\quad\in S_n,
\end{equation}
we shall prove in a moment that (with obvious notation)
\begin{equation}\label{EqDistZr'}
\Delta_{{\cal F}^d_{S,per}(s,R)}(Z_r,\bar Z_r)\lesssim
n^{-1/2}n_0+\sigma^{-1}Rn_0^{1/2-s/d},
\end{equation}
but then the assertion follows: Observing $\bar Z_r$ is equivalent
to observing
\[ dY(x)=P_nf(x)+\sigma n^{-1/2}dB(x),\quad x\in [0,1]^d,\]
which has a total variation distance to the Gaussian shift $\G_n^d$
of order $\sigma^{-1}n^{1/2}\norm{f-P_nf}_{L^2}\lesssim
\sigma^{-1}n^{1/2-s/d}\norm{f}_{H^s}$. Using the triangle inequality
for the Le Cam distance between the intermediate experiments, we
arrive at the bound for $\Delta_{{\cal
F}^d_{S,per}(s,R)}(\E^d_{n,r},\G^d_n)$.

To obtain \eqref{EqDistZr'}, we take expectations over the design
and split
\[E[\norm{{\cal L}(Z_r)-{\cal L}(Z_r')}_{TV}^2{\bf
1}_{\Omega_{n_0}^n}]\lesssim I+II+III
\]
with the terms
\begin{align*}
I&:=n\sigma^{-2}E\big[\norm{(P_{n_0}^n-P_{n_0})f}_{L^2}^2{\bf
1}_{\Omega_{n_0}^n}\big] \text{ (difference in
mean on $S_{n_0}$)},\\
II&:=E[\norm{(\Pi_n|_{S_{n_0}})^{-1}-\Id_{S_{n_0}}}_{HS}^2{\bf
1}_{\Omega_{n_0}^n}] \text{ (heteroskedasticity on
$S_{n_0}$)},\\
III&:=n\sigma^{-2}E\big[\norm{(T^{-1}(P_n^n-P_{n_0}^n)-(P_n-P_{n_0}))f}_{L^2}^2{\bf
1}_{\Omega_{n_0}^n}\big] \text{ (difference in mean on
$S_{n_0}^{\perp_{L^2}}$)}.
\end{align*}

{\noindent\bf Term I.} Using the projection properties, we obtain on
$\Omega_{n_0}^n$:
\[ \norm{(P_{n_0}^n-P_{n_0})f}_{L^2}^2=\norm{P_{n_0}^n(\Id-P_{n_0})f}_{L^2}^2
\le 4\norm{P_{n_0}^n(\Id-P_{n_0})f}_{n}^2.
\]
Because of
$E[\scapro{\phi_k}{\phi_j^n}_n\overline{\scapro{\phi_{k'}}{\phi_j^n}_n}]=0$
for $k\not=k'$, $k,k'>j$ by Proposition \ref{PropE=0} below, an
expansion in the basis $(\phi_j^n)$ yields
\begin{align*}
E\big[\norm{P_{n_0}^n(\Id-P_{n_0})f}_{n}^2\big]&=
\sum_{j=1}^{n_0}\sum_{k=n_0+1}^\infty\abs{\scapro{f}{\phi_k}_{L^2}}^2
E[\abs{\scapro{\phi_k}{\phi_j^n}_n}^2]\\
&=\sum_{k=n_0+1}^\infty\abs{\scapro{f}{\phi_k}_{L^2}}^2E[\norm{P^n_{n_0}\phi_k}_n^2].
\end{align*}
Proposition \ref{PropProjEst} below yields
$E[\norm{P^n_{n_0}\phi_k}_n^2]\lesssim k/n$ and hence
\[ I\lesssim
\sigma^{-2}\sum_{k=n_0+1}^\infty\abs{\scapro{f}{\phi_k}_{L^2}}^2k\lesssim
\sigma^{-2}n_0^{1-2s/d}\norm{f}_{H^s}^2.
\]

{\noindent\bf Term II.} Using
$\norm{(\Pi_n|_{S_{n_0}})^{-1}}_{L^2\to L^2}\le 4$ on
$\Omega_{n_0}^n$, we find:
\begin{align*}
E[\norm{(\Pi_n|_{S_{n_0}})^{-1}-\Id_{S_{n_0}}}_{HS}^2{\bf
1}_{\Omega_{n_0}^n}]&\le E[\norm{(\Pi_n|_{S_{n_0}})^{-1}}_{L^2\to
L^2}\norm{\Pi_n|_{S_{n_0}}-\Id_{S_{n_0}}}_{HS}^2
{\bf 1}_{\Omega_{n_0}^n}]\\
&\le 4E[\norm{\Pi_n|_{S_{n_0}}-\Id_{S_{n_0}}}_{HS}^2]\\
&=4\sum_{j,j'=1}^{n_0}
E[\abs{\scapro{\phi_j}{\phi_{j'}}_n-\delta_{j,j'}}^2]\\
&\le 4n^{-1}\sum_{j,j'=1}^{n_0} \int
\abs{\phi_j}^2\abs{\phi_{j'}}^2.
\end{align*}
For the Fourier basis we obtain $II\le 4n^{-1}n_0^2$.\\

{\noindent\bf Term III.} Let us write $f=f_0+f_1+f_2$ with
$f_0=P_{n_0}f$, $f_1=(P_n-P_{n_0})f$, $f_2=(\Id-P_n)f$. Then the
projection properties imply
\begin{align*}
&E\big[\norm{(T^{-1}(P_n^n-P_{n_0}^n)-(P_n-P_{n_0}))f}_{L^2}^2{\bf 1}_{\Omega_{n_0}^n}\big]\\
&=E\big[\norm{T^{-1}f_1+T^{-1}P_n^nf_2-T^{-1}P_{n_0}^n(f_1+f_2)-f_1}_{L^2}^2{\bf 1}_{\Omega_{n_0}^n}\big]\\
&\le
3E\big[\norm{(T^{-1}-\Id)f_1}_{L^2}^2+\norm{(P_n^n-P_{n_0}^n)f_2}_n^2
+\norm{P^n_{n_0}f_1}_n^2{\bf 1}_{\Omega_{n_0}^n}\big]\\
&\le 3E\big[\norm{f_1}_{n}^2+\norm{f_1}_{L^2}^2-
2\Re(\scapro{T^{-1}f_1}{f_1}_{L^2})\big]+3E\big[\norm{f_2}_n^2\big]
+3E\big[\norm{P^n_{n_0}f_1}_n^2{\bf 1}_{\Omega_{n_0}^n}\big]\\
&=6E\big[\Re(\scapro{f_1-T^{-1}f_1}{f_1}_{L^2})\big]+3\norm{f_2}_{L^2}^2
+3E\big[\norm{P_{n_0}^nf_1}_{n}^2{\bf 1}_{\Omega_{n_0}^n}\big]\\
&=:III_1+III_2+III_3.
\end{align*}

The term $III_2$ is easily bounded by $\norm{f_2}_{L^2}^2\lesssim
n^{-2s/d}\norm{f}_{H^s}^2$. As in the estimate for term I, we obtain
$III_3\lesssim n^{-1}n_0^{1-2s/d}\norm{f}_{H^s}^2$. For $III_1$ we
use $E[\scapro{T^{-1}\phi_j}{\phi_k}_{L^2}]=0$, $j\not=k$, by
Proposition \ref{PropE=0} below to conclude
\[E\big[\Re(\scapro{f_1-T^{-1}f_1}{f_1}_{L^2})\big]=\sum_{j=n_0+1}^n\abs{\scapro{f}{\phi_j}_{L^2}}^2
E\big[\scapro{(\Id-T^{-1})\phi_j}{\phi_j}_{L^2}].
\]
Because of $\norm{\phi_j}_n=1$ for the Fourier basis we find
\[\scapro{T^{-1}\phi_j}{\phi_j}_{L^2}
=\scapro{\norm{\phi_j-P_{j-1}^n\phi_j}_n\phi_j^n+P_{j-1}^n\phi_j}{\phi_j^n}_n
=\norm{\phi_j-P_{j-1}^n\phi_j}_n\ge 1-\norm{P_{j-1}^n\phi_j}_n^2.
\]
By Proposition \ref{PropProjEst} below, the bound
\[
E\big[\Re(\scapro{f_1-T^{-1}f_1}{f_1}_{L^2})\big] \le
\sum_{j=n_0+1}^{n}\abs{\scapro{f}{\phi_j}_{L^2}}^2
E[\norm{P_{j-1}^n\phi_j}_n^2] \lesssim
\sum_{j=n_0+1}^{n}\frac{j}{n}\abs{\scapro{f}{\phi_j}_{L^2}}^2
\]
follows, which is of order $n^{-1}n_0^{1-2s/d}\norm{f}_{H^s}^2$.
Putting the estimates together, we have shown
\[ III\lesssim \sigma^{-2}\big(
n_0^{1-2s/d}\norm{f}_{H^s}^2+n^{1-2s/d}\norm{f}_{H^s}^2+n_0^{1-2s/d}\norm{f}_{H^s}^2\big)\lesssim
\sigma^{-2}n_0^{1-2s/d}\norm{f}_{H^s}^2
\]
and in sum $I+II+III\lesssim \sigma^{-2}n_0^{1-2s/d}R^2+n^{-1}n_0^2$
uniformly over $f\in {\cal F}^d_{S,per}(s,R)$, which gives the
asserted bound \eqref{EqDistZr'}.
\end{proof}

\subsection{Technical results}

We gather results on fine properties of the Fourier basis $(\phi_j)$
and its generated approximation spaces $S_n$. The setting is as in
the proof of Theorem \ref{ThmRanDes}. For the value of the next
proposition notice that
$\scapro{\phi_{k'}}{\phi_k^n}_n=\scapro{T^{-1}\phi_{k'}}{\phi_k}_{L^2}$.

\begin{proposition}\label{PropE=0}
We have for indices  $k'',k'>k\ge 1$, $k''\not=k'$:
\[E[\scapro{\phi_{k'}}{\phi_k^n}_n]=0\qquad \text{and}\qquad
E[\scapro{\phi_{k'}}{\phi_k^n}_n\overline{\scapro{\phi_{k''}}{\phi_k^n}}_n]=0.
\]
\end{proposition}

\begin{proof}
Since the randomness enters via $P_{k-1}^n$ in a very intricate way,
we use a symmetry argument. Specify $X_i:=(Y_i+\theta) \mod 1$,
$i=1,\ldots,n$, with $Y_i\sim U([0,1]^d)$, $\theta\sim U([0,1]^d)$
all independent such that $X_i\sim U([0,1]^d)$ i.i.d. Working
conditionally on $\theta$, we shall keep track on the dependence on
$\theta$ using brackets. We claim that for $k'>k$ it holds that
\begin{equation}\label{EqInv}
\scapro{\phi_{k'}}{\phi_k^n}_n[\theta]=e^{2\pi
\iota\scapro{\ell(k')-\ell(k)}{\theta}}\scapro{\phi_{k'}}{\phi_k^n}_n[0],
\end{equation}
which entails the result due to
\[\int_{[0,1]^d} e^{2\pi \iota
\scapro{\ell(k')-\ell(k)}{\theta}}d\theta=0\text{ and }
\int_{[0,1]^d} e^{2\pi \iota (\scapro{\ell(k')-\ell(k)}{\theta}-
\scapro{\ell(k'')-\ell(k)}{\theta})}d\theta=0.
\]
For $m\in\Z^d$ put
\[ A_m[\theta]:=\frac1n\sum_{j=1}^n e^{2\pi \iota\scapro{m}{X_j[\theta]}}=\frac1n\sum_{j=1}^n e^{2\pi
\iota\scapro{m}{Y_j+\theta}}=e^{2\pi \iota\scapro{m}{\theta}}A_m[0].
\]
The proof of \eqref{EqInv} will be performed by induction from
$\kappa<k$ to $k$, considering tupels $(\kappa',\kappa)$,
$\kappa'>\kappa$, and $(k',k)$, $k'>k$. Since $\ell(1)=0$ and
$\phi_1^n=\phi_1=1$, we have for $k'>1$ and $k=1$
\[\scapro{\phi_{k'}}{\phi_k^n}_n[\theta]=\frac1n\sum_{j=1}^n
e^{2\pi \iota \scapro{\ell(k')}{Y_j+\theta}}=e^{2\pi \iota
\scapro{\ell(k')-\ell(1)}{\theta}}\scapro{\phi_{k'}}{\phi_k^n}_n[0].
\]
Writing $c_k:=\norm{\phi_k-P_{k-1}^n\phi_k}_n^{-1}$, the induction
hypothesis implies
\[c_k^{-2}[\theta]=1-\sum_{j=1}^{k-1}\abs{\scapro{\phi_k}{\phi_j^n}}^2[\theta]=c_k^{-2}[0]
\]
and furthermore
\begin{align*}
\scapro{\phi_{k'}}{\phi_k^n}_n[\theta] &=
\scapro{\phi_{k'}}{c_k(\phi_k-\textstyle\sum_{r=1}^{k-1}\scapro{\phi_k}{\phi_r^n}\phi_r^n)}[\theta]\\
&=c_k\Big(\scapro{\phi_{k'}}{\phi_k}_n-\sum_{r=1}^{k-1}\scapro{\phi_{k'}}{\phi_r^n}_n
\overline{\scapro{\phi_k}{\phi_r^n}_n}\Big)[\theta]\\
&=c_k[\theta]\Big(A_{\ell(k')-\ell(k)}[\theta]
-\sum_{r=1}^{k-1}e^{2\pi \iota \scapro{\ell(k')-\ell(k)}{\theta}}
\scapro{\phi_{k'}}{\phi_r^n}_n[0]
\overline{\scapro{\phi_k}{\phi_r^n}_n[0]}\Big)\\
&=e^{2\pi \iota \scapro{\ell(k')-\ell(k)}{\theta}}
\scapro{\phi_{k'}}{\phi_k^n}_n[0],
\end{align*}
which proves the induction step and gives \eqref{EqInv}.
\end{proof}

\begin{proposition}\label{PropFourierDiscr}
Suppose $g=\sum_{\abs{\ell}_{\ell^2}\le L}\gamma_\ell e^{2\pi\iota
\scapro{\ell}{\cdot}}$ is a $d$-dimensional trigonometric polynomial
of degree $L$. Let $\Delta\in (0,L^{-1}]$ with $1/\Delta\in\N$ be
given and define the cubes $C_m:=\prod_{i=1}^d
[(m_i-1)\Delta,m_i\Delta)$. Then:
\[ \Delta^d\sum_{m\in\{1,\ldots,\Delta^{-1}\}^d}\sup_{x_m\in
C_m}\babs{\abs{g(x_m)}^2-\abs{g(m\Delta)}^2}\le
\norm{g}_{L^2}^2(e^{2d\Delta L}-1).
\]
\end{proposition}

\begin{proof}
We need multi-indices $\alpha,\beta\in\N_0^d$ with
$\alpha!:=\alpha_1!\cdots\alpha_d!$, $x^\alpha:=x_1^{\alpha_1}\cdots
x_d^{\alpha_d}$,
$\binom{\alpha}{\beta}:=\frac{\alpha!}{\beta!(\alpha-\beta)!}$ and
differential operators
$D^\alpha:=\frac{\partial^{\alpha_1}}{\partial x_1^{\alpha_1}}\cdots
\frac{\partial^{\alpha_d}}{\partial x_d^{\alpha_1}}$. Since
$\abs{g}^2$ is real-analytic, a power series expansion gives for any
$x_m\in C_m$:
\begin{align*}
\abs{\abs{g}^2(x_m)-\abs{g}^2(m\Delta)}&=\babs{\sum_{\alpha\in\N_0^d,\alpha\not=0}
D^\alpha\abs{g}^2(m\Delta)\frac{(x_m-m\Delta)^\alpha}{\alpha!}}\\
&\le \sum_{\alpha\in\N_0^d,\alpha\not=0}
\frac{\Delta^{\abs{\alpha}_{\ell^1}}}{\alpha!}
\sum_{\beta\in\N_0^d,\beta\le \alpha} \binom{\alpha}{\beta}
\abs{D^\beta g(m\Delta)} \abs{D^{\alpha-\beta}\bar{g}(m\Delta)}.
\end{align*}
Together with $g$ any derivative is again a trigonometric polynomial
of degree $L$ and by the isometry \eqref{EqFourierScaPro} and
Bernstein's inequality, cf. \cite[p. 32]{Meyer}, we obtain
\[ \Delta^{d}\sum_{m\in\{1,\ldots,\Delta^{-1}\}^d} \abs{D^\alpha
g}^2(m\Delta)=\norm{D^\alpha g}_{L^2}^2\le
L^{2\abs{\alpha}_{\ell^1}}\norm{g}_{L^2}^2.
\]
This implies by the Cauchy-Schwarz inequality
\begin{align*}
&\Delta^d\sum_{m\in\{1,\ldots,\Delta^{-1}\}^d}\sup_{x_m\in
C_m}\abs{\abs{g}^2(x_m)-\abs{g}^2(m\Delta)}\\
&\le \Delta^d\sum_{\alpha\in\N^d,\alpha\not=0}
\frac{\Delta^{\abs{\alpha}_{\ell^1}}}{\alpha!}
\sum_{\beta\in\N^d,\beta\le \alpha} \binom{\alpha}{\beta}
\Big(\sum_{m}\abs{D^\beta g(m\Delta)}^2\Big)^{1/2}
\Big(\sum_{m}\abs{D^{\alpha-\beta}\bar{g}(m\Delta)}^2\Big)^{1/2}\\
&\le \norm{g}_{L^2}^2 \sum_{\alpha\in\N^d,\alpha\not=0}
\frac{\Delta^{\abs{\alpha}_{\ell^1}}}{\alpha!}
\sum_{\beta\in\N^d,\beta\le\alpha} \binom{\alpha}{\beta}
L^{\abs{\beta}_{\ell^1}}L^{\abs{\alpha-\beta}_{\ell^1}}\\
&= \norm{g}_{L^2}^2 \sum_{\alpha\in\N^d,\alpha\not=0}
\frac{(2\Delta L)^{\abs{\alpha}_{\ell^1}}}{\alpha!}\\
&= \norm{g}_{L^2}^2(e^{2d\Delta L}-1).
\end{align*}

\end{proof}

\begin{lemma}\label{LemMaxMult}
Let $Y\in\R^r$ follow the multinomial distribution with parameters
$n$ and $p_1=\cdots=p_r=1/r$. Then for $n\to\infty$ and $r=r(n)$
with $r\log(r)/n\to 0$
\[ \forall\,C>0:\:\limsup_{n\to\infty}\tfrac14 r(n)^{C^2/4-1}P\Big(\max_{1\le i\le
r(n)}\abs{Y_i-n/r(n)}>C\sqrt{n\log(r(n))/r(n)}\Big)\le 1.
\]
\end{lemma}

\begin{proof}
If $X_1,\ldots,X_r$ are independently Poisson($n/r$)-distributed,
then it is well known that the law of $(X_1,\ldots,X_r)$ given
$\sum_{i=1}^rX_i=n$ is multinomial with parameters $n$ and
$p_1=\cdots=p_r=1/r$. Set $A_{nr}:=C\sqrt{n\log(r)/r}$. Since
\[\textstyle k\mapsto P\Big(\max_{1\le i\le
r}X_i-n/r>A_{nr}\,\Big|\,\sum_{i=1}^rX_i=k\Big)
\]
is obviously increasing in $k\in\N$, we obtain
\[ {\textstyle P\Big(\max_{1\le i\le
r}X_i-n/r>A_{nr}\,\Big|\,\sum_{i=1}^rX_i=n\Big)}\le
\frac{P\big(\max_{1\le i\le
r}X_i-n/r>A_{nr}\big)}{P\big(\sum_{i=1}^rX_i\ge n\big)}.
\]
As $\sum_{i=1}^rX_i$ is Poisson($n$)-distributed,
$\lim_{n\to\infty}P(\sum_{i=1}^rX_i\ge n)=1/2$ holds, whence
\begin{equation}\label{EqXY}
\textstyle\limsup_{n\to\infty}\Big(P\big(\max_{1\le i\le
r}Y_i-n/r>A_{nr}\big)-2P\big(\max_{1\le i\le
r}X_i-n/r>A_{nr}\big)\Big)\le 0.
\end{equation}
By the exponential moment estimate
$E[e^{a(X_i-n/r)}]=e^{n(e^a-a-1)/r}\le e^{3na^2/4r}$ for
$a:=rA_{nr}/n\to 0$ and $n$ large, the generalized Markov inequality
yields
\[\textstyle P(\max_{1\le i\le r}X_i-n/r>A_{nr})\le rP(X_i-n/r>A_{nr})\le
re^{3na^2/4r-aA_{nr}}=r^{1-C^2/4}.
\]
By use of \eqref{EqXY} and a completely symmetric argument for
$P(\max_{1\le i\le r}(n/r-X_i)>A_{nr})$, the result follows.
\end{proof}

\begin{proposition}\label{PropOmegaj}
For $j=j(n)$ such that $j\log(j)=o(n)$ and the event $\Omega_j^n$ in
\eqref{EqOmegaj} we have
$\lim_{n\to\infty}n^pP((\Omega_{j(n)}^n)^\complement)=0$ for any
power $p>0$.
\end{proposition}

\begin{proof}
From Proposition \ref{PropFourierDiscr} we derive with $\Delta\le
L:=\abs{\ell(j)}_{\ell^2}$, $1/\Delta\in\N$, the cubes
$C_m:=\prod_{i=1}^d [(m_i-1)\Delta,m_i\Delta)$ and the occupations
$N_m:=\#\{i:X_i\in C_m\}$:
\begin{align*}
&\babs{\norm{g}_{L^2}^2-\frac1n\sum_{i=1}^n\abs{g(X_i)}^2}\\
&= \frac1n\babs{\sum_{m\in\{1,\ldots,\Delta^{-1}\}^d}\Big(\Delta^d
n\abs{g(m\Delta)}^2-\sum_{i:X_i\in C_m} \abs{g(X_i)}^2\Big)}\\
&\le \frac1n\sum_{m\in\{1,\ldots,\Delta^{-1}\}^d}\Big(\abs{\Delta^d
n-N_m}\abs{g(m\Delta)}^2+N_m\sup_{x_m\in C_m} \babs{\abs{g(m\Delta)}^2-\abs{g(x_m)}^2}\Big)\\
&\le \frac{\norm{g}_{L^2}^2}{\Delta^{d}n}
\max_{m\in\{1,\ldots,\Delta^{-1}\}^d} \Big(\abs{\Delta^d
n-N_m}+N_m(e^{2d\Delta L}-1)\Big)\\
&\le \norm{g}_{L^2}^2 \Big( e^{2d\Delta L}
\max_{m\in\{1,\ldots,\Delta^{-1}\}^d}
\abs{1 -N_m/n\Delta^d} + (e^{2d\Delta L}-1) \Big)\\
\end{align*}
By Lemma \ref{LemMaxMult} $\max_m\abs{1-N_m/n\Delta^d}^2\ge C
(n\Delta^d)^{-1}\log(1/\Delta)$ has probability tending to zero with
any given polynomial rate when choosing $C$ sufficiently large.
Since $L^d\log(L)\lesssim j\log(j)=o(n)$, we can choose
$\Delta=o(L^{-1})$ such that still
$\Delta^{-d}\log^2(1/\Delta)=o(n)$ holds. This gives
\[\abs{\norm{g}_{L^2}^2-\norm{g}_n^2}\le \Big( C e^{2d\Delta L}
(n\Delta^d)^{-1}\log(1/\Delta)+ (e^{2d\Delta L}-1)
\Big)\norm{g}_{L^2}^2\le \tfrac34\norm{g}_{L^2}^2
\]
for large $n$ with probability larger than $1-n^{-p}$.
\end{proof}

\begin{proposition} \label{PropProjEst}
For $j\in\N$ with $j\log(j)=o(n)$ we have
\[ E[\norm{P^n_{j-1}\phi_j}_n^2]\lesssim j/n.\]
\end{proposition}

\begin{proof}
By construction, $\norm{P^n_{j-1}\phi_j}_n^2\le \norm{\phi_j}_n^2=1$
holds so that by Proposition \ref{PropOmegaj} it suffices to find
the bound for the expectation on the event $\Omega_j^n$.


Setting $A_m:=\frac1n\sum_{k=1}^m \exp(2\pi \iota\scapro{m}{X_k})$,
$m\in\Z^d$, we use Parseval's identity and $E[\abs{A_m}^2]=1/n$ for
$m\not=0$ to obtain
\begin{align*}
E[\norm{P^n_{j-1}\phi_j}_n^2{\bf 1}_{\Omega_j^n}] &=
E\Big[\sup_{g\in
V_{j-1}}\frac{\abs{\scapro{\phi_j}{g}_n}^2}{\norm{g}_n^2}{\bf 1}_{\Omega_j^n}\Big]\\
&\le E\Big[\sup_{\norm{(c_r)}_{\ell^2}=1} \babs{\frac1n\sum_{k=1}^n
\sum_{r=1}^{j-1}\bar c_re^{2\pi \iota
\scapro{\ell(j)-\ell(r)}{X_k}}}^2 \sup_{g\in
V_{j-1}}\frac{\norm{g}_{L^2}^2}{\norm{g}_n^2}{\bf 1}_{\Omega_j^n}\Big]\\
&\le 4E\Big[\sum_{r=1}^{j-1}\abs{A_{\ell(j)-\ell(r)}}^2
\Big]=\frac{4(j-1)}{n}.
\end{align*}
This gives the result.
\end{proof}

\bibliographystyle{economet}

\bibliography{mybib}

\end{document}